\documentclass[12pt,oneside]{amsart}
\usepackage{latexsym}
\usepackage{amsmath,amssymb,amsthm}
\usepackage{amscd}
\usepackage[dvips]{graphics}

\newcommand{\re}{\mathbb{R}}
\newcommand{\co}{\mathbb{C}}
\newcommand{\cp}{\mathbb{CP}}

\newcommand{\Sl}[1]{\mathbf{SL}(#1,\mathbb{R})}
\newcommand{\gl}[1]{\mathbf{GL}(#1,\mathbb{R})}
\newcommand{\na}{\nabla}
\newcommand{\sfrac}[2]{{\textstyle \frac{#1}{#2}}}
\newcommand{\D}{\displaystyle}

\newtheorem{prop}{Proposition} 
\newtheorem{cor}[prop]{Corollary}
\newtheorem{thm}{Theorem}

\newtheorem{lem}[prop]{Lemma}

\theoremstyle{remark}
\newtheorem*{rem}{Remark}
\newtheorem*{ack}{Acknowledgements}

\begin{document}
\title{Singular Semi-Flat Calabi-Yau Metrics on $S^2$}
\author{John C. Loftin}
 \maketitle


\section{Introduction}

This paper is motivated by recent work of Gross and Wilson
\cite{gross-wilson00}, in which they construct degenerate limits
of families of K3 surfaces equipped with Calabi-Yau metrics
(Ricci-flat K\"ahler metrics).  Upon proper rescaling, the metric
limit of such a family is a two-dimensional sphere equipped with a
Riemannian metric with prescribed singularities at 24 points. Away
from singularities, this limit metric is an affine K\"ahler
metric.  In other words, there are natural affine flat coordinates
$(\alpha^j)$ and a local potential function $\phi$ so that the
metric is given by
\begin{equation}\label{hess-met}
\frac{\partial^2\phi}{\partial \alpha^j \partial \alpha^k} \,
d\alpha^j\, d \alpha^k.
\end{equation}
 Moreover, the potential $\phi$ satisfies the Monge-Amp\`ere
 equation
\begin{equation}\label{m-a}
\det \frac{\partial^2\phi}{\partial \alpha^j \partial \alpha^k} =
1.
\end{equation}
In this case, the metric is naturally a real slice of a Calabi-Yau
metric.  We refer to such a metric as a \emph{semi-flat Calabi-Yau
metric}.

Such singular semi-flat Calabi-Yau metrics on surfaces were first
constructed by Greene-Shapere-Vafa-Yau \cite{greene-svy90}. We
construct many examples of such metrics.  Our main theorem is
this:
\begin{thm}
Given any holomorphic cubic differential $U$ on $\cp^1$ which has
poles of order 1 at a finite number of points $p_j$ and is not
identically zero, there exists an affine flat structure and a
semi-flat Calabi-Yau metric on $\cp^1 \setminus \{p_j\}$.  The
singularities of the affine flat structure and metric at the $p_j$
are asymptotically the same as those in \cite{gross-wilson00}.
\end{thm}
The details of the nature of the singularities are given in
Theorem \ref{big-theorem} below.

A semi-flat Calabi-Yau metric can naturally be seen as a real
slice of a Calabi-Yau metric.  For a K\"ahler metric given by a
potential function $\Phi$, the equation for the metric to be Ricci
flat is
$$ \partial \bar \partial \log \det \Phi_{j\bar k} = 0.$$
A convex function $\phi$ on a domain $\Omega \subset \re^n$ can be
extended to be constant along the imaginary fibers of the tube
domain $\Omega + i\re^n\subset \co^n$.  Then up to a constant
factor, the K\"ahler metric $\phi_{j\bar k} \,dz^j\overline{dz^k}$
extends the affine K\"ahler metric $\phi_{jk}\,d\alpha^j
d\alpha^k$ on $\Omega$.  If in addition $\phi$ satisfies $\det
\phi_{jk}=1$, then $\phi_{j\bar k} \,dz^j\overline{dz^k}$ is
Calabi-Yau. This K\"ahler metric is commonly known as semi-flat
since it is flat along the imaginary fibers. We then extend the
terminology a small bit to call the real metric $\phi_{jk} \,
d\alpha^jd\alpha^k$ semi-flat Calabi-Yau in this case.  The
conjecture of Strominger-Yau-Zaslow \cite{syz} (as explicated by
Hitchin, Gross, Leung and others) implies that in this degenerate
limit, mirror symmetry reduces to the Legendre transform of the
affine K\"ahler potential function $\phi$.

The graph in $\re^{n+1}$ of a convex function $\phi$ satisfying
(\ref{m-a}) is a hypersurface classically studied in affine
differential geometry, a parabolic affine sphere.  The study of
such surfaces dates back to \c{T}i\c{t}eica and Blaschke.  A
parabolic affine sphere whose metric (\ref{hess-met}) is complete
must be the graph of a quadratic polynomial (and thus the metric
is flat).  This was proved by J\"orgens \cite{jorgens54} in the
case $n=2$ and by Calabi \cite{calabi58} in higher dimensions.
Cheng and Yau \cite{cheng-yau82} studied affine K\"ahler manifolds
and produced semi-flat Calabi-Yau metrics on many compact
manifolds (on any compact affine K\"ahler manifold which admits a
volume form covariantly constant with respect to the the canonical
affine flat connection).  We note affine K\"ahler metrics are also
studied in the works of Shima e.g.\ \cite{shima81}, where they are
called Hessian metrics.  Cheng and Yau observed that Calabi's
estimates imply all semi-flat Calabi-Yau metrics on a compact
manifold are flat. Thus in the present work it is important that
the metric we produce is not complete near the singularities: if
it were, then Calabi's theorem would imply that it is flat. We
should also mention that the Bernstein problem has for parabolic
affine spheres has also been solved by J\"orgens ($n=2$), Calabi
($n=3,4,5$), Pogorelov \cite{pogorelov72}, and Cheng-Yau
\cite{cheng-yau86}.

Given an affine K\"ahler metric with local potential function
$\phi$, the affine K\"ahler metric transforms as a tensor under
affine flat coordinate changes. A manifold $M$ built of coordinate
charts in $\re^n$ with gluing maps locally constant maps in
$\gl{n}\ltimes\re^n$ is called an affine flat manifold. This
system of canonical affine flat coordinates is equivalent to the
existence of a torsion-free flat connection $\na$ on the tangent
bundle (the coordinate vector fields $\partial/\partial \alpha^j$
are parallel with respect to this connection).  In our case,
equation (\ref{m-a}) demands a little more structure on the
manifold $M$.  The 1 on the right side of (\ref{m-a}) is parallel
under $\na$ and is the square of a volume form. Thus in general
there is a parallel density on $M$, and if $M$ is oriented (as are
the examples in this paper) $M$ admits a parallel volume form.
Then the gluing maps are elements of $\Sl{n}\ltimes\re^n$. We can
think of $\na$ as an affine connection (a connection on a
principal bundle modelled on the affine group
$\gl{n}\ltimes\re^n$). Thus we have a natural holonomy
representation $\gamma$ on an affine flat manifold $M$ with
parallel volume form: $$\gamma\!: \pi_1M \to \Sl{n}\ltimes\re^n.$$
In the present work, we calculate the conjugacy class of the
holonomy around each singular point.

In dimension $n=2$, there are two special techniques for analyzing
parabolic affine spheres.  Both use the conformal structure
induced by the affine metric (\ref{hess-met}). The first involves
a semilinear equation used by Simon and Wang \cite{simon-wang93}
to perform a conformal change to find the affine metric.  (We call
this semilinear equation \c{T}i\c{t}eica's equation.) Then Simon
and Wang introduce an initial-value problem (a developing map) to
deduce the local structure of the parabolic affine sphere. The
second technique, which goes back to Blaschke, is a representation
of the parabolic affine sphere in terms of two holomorphic
functions. We recall both these theories, as well as some general
facts about parabolic affine spheres and affine flat coordinates,
in Section \ref{parabolic} below.

In Section \ref{solve-u}, we introduce a model metric solution to
\c{T}i\c{t}eica's equation near each first-order pole of a cubic
differential $U$ on $\cp^1$.  This model comes from Gross-Wilson
\cite{gross-wilson00}.  Then we proceed to perturb the model
metric by a conformal factor $e^u$ and get asymptotic bounds on
$u$ near each singularity.

Then in Section \ref{affine-flat}, using the bounds on $u$, we
analyze the affine flat structure induced on
$M=\cp^1\setminus\{p_j\}$ induced by the parabolic affine sphere
we've constructed. In particular, we use Simon-Wang's developing
map and techniques of ODEs to calculate the holonomy and other
natural invariants of affine flat structure.  (This basic plan of
first solving for a conformal factor and then applying a
developing map and ODE techniques to characterize the relevant
geometric structure near a singular point on a surface was first
carried out in \cite{loftin02c}, where asymptotics for singular
convex real projective structures were investigated using
hyperbolic affine spheres.) Finally we use Blaschke's holomorphic
characterization to find precise asymptotics of the metric and the
affine flat structure.

Then in Section \ref{mirror}, we recall Leung's picture of mirror
symmetry without correction terms \cite{leung00,leung02}, and we
write down the manner in which mirror symmetry should work in this
degenerate limit via the Legendre transform.

We note that we are not able to reproduce one relevant feature
required by mirror symmetry, the integrality of the affine
holonomy. In the picture of Strominger-Yau-Zaslow, the total space
of a Calabi-Yau manifold is formed by a fibration of special
Lagrangian tori over a singular affine flat manifold $B$. (In the
present work, $B$ is $S^2=\cp^1$ with singularities at the $p_j$.)
These tori are naturally quotients of the imaginary fibers of the
tube domain construction above.  In order for such a quotient to
make sense globally on an affine flat manifold, the linear part of
holonomy should be integral (so that there is a lattice in the
tangent bundle preserved by $\na$). In other words, the holonomy
representation should be conjugate to one in the group
$\mathbf{SL}(2,{\mathbb Z})\ltimes \re^2$.  In this paper we
determine the holonomy only near the singular points, and it seems
the methods of this paper are insufficient to determine such a
global integrality condition.  On the other hand, there are
combinatorial constructions of integral affine manifolds with
singularities due to Haase-Zharkov \cite{haase-zharkov02} and
Gross-Siebert \cite{gross-siebert03a,gross-siebert03b}, which
discuss mirror symmetry from combinatorial and algebro-geometric
points of view.  Haase-Zharkov \cite{haase-zharkov03} also
construct affine K\"ahler metrics on their examples, but these do
not satisfy the Monge-Amp\`ere equation.  Recently Zharkov
\cite{zharkov03} has conjectured a detailed picture of how
Calabi-Yau metrics degenerate to semi-flat Calabi-Yau metrics.

\begin{ack} The author would like to thank many people for
valuable discussions. Among them are Conan Leung and Jacob Sturm.
The author would also like to thank Andreea Nicoara for pointing
out the correct spelling of \c{T}i\c{t}eica, and Robert Bryant for
his remark relating semi-flat Calabi-Yau metrics to hyperbolic
metrics in dimension 2. Finally, the author would like to thank
Eric Zaslow and S.T. Yau for many useful discussions.
\end{ack}

\section{Parabolic affine spheres} \label{parabolic}

\subsection{Affine flat coordinates} \label{aff-flat-coord}

Given any locally strictly convex immersed hypersurface $H\subset
\re^{n+1}$, there is a natural transversal vector field
$\xi=\xi_H$ which is invariant under affine volume-preserving
automorphisms of $\re^{n+1}$ (elements of the group
$\Sl{n+1}\ltimes\re^{n+1}$).  $\xi$ is called the \emph{affine
normal} to $H$. In other words, if $\Psi
\in\Sl{n+1}\ltimes\re^{n+1}$ and $p\in H$, then $$ \Psi_* \xi_H(p)
= \xi_{\Psi(H)}(\Psi(p)).$$ $H$ is called a \emph{parabolic affine
sphere} if $\xi$ is a constant vector. It is standard to choose
coordinates on $\re^{n+1}$ so that $\xi =(0,\dots,0,1)$ in this
case.

In these coordinates $H$ can be locally represented as the graph
of a strictly convex function $\phi$---so that
$H=\{(\alpha,\phi(\alpha))\}$ for $\alpha$ in a domain in $\re^n$.
The condition that $H$ be a parabolic affine sphere is then the
real Monge-Amp\`ere equation
$$ \det \frac{\partial^2\phi}{\partial \alpha^j \partial \alpha^k} = 1.$$

The Legendre transform is quite natural in this context. Recall
that if $\beta_j=\partial \phi / \partial \alpha^j$, then the
Legendre transform $\chi$ of $\phi$ is given by
$$\chi + \phi = \beta_j \alpha^j.$$
(We use the usual summation convention.) We primarily think of
$\chi$ as a function of the $\beta_j$, which are coordinates on
the dual vector space $\re_n$ of $\re^n$.   The graph
$(\beta,\chi)\subset \re_{n+1}$ is again a parabolic affine
sphere.

The natural group action in this setting consists of those
elements $\Psi\in\Sl{n+1}\ltimes\re^{n+1}$ which preserve $\xi$ in
the sense that $\Psi_*\xi=\xi$.  In other words, we are interested
in the group $G$ of transformations of the form
$$(\tilde \alpha \,\, \tilde \gamma) =
\Psi (\alpha \,\,\gamma) = (\alpha \,\, \gamma) \left(
     \begin{array}{cc} A & c \\ 0 & 1 \end{array} \right)
     + (b \,\, d),$$
where $\alpha$ and $b$ are row vectors in $\re^n$, $c$ is a column
vector, $\gamma$ and $d \in \re$, and $A\in \Sl{n}$.

The Legendre transform is natural with respect to these
coordinates.  In terms of the $(\tilde \alpha \,\, \tilde \gamma)$
coordinates, $H$ is the graph of a function $\tilde \phi$ so that
$\tilde\gamma = \tilde \phi(\tilde\alpha)$.  Form the Legendre
transform: $\tilde\beta_j = \partial \tilde\phi /
\partial \tilde\alpha^j$ and $\tilde\chi+\tilde \phi = \tilde
\beta_j\tilde \alpha^j$.  Then it is straightforward to verify
that
$$ \left(\begin{array}{c} \tilde \beta \\ \tilde \chi
\end{array}\right)
= \left( \begin{array}{cc} A^{-1} & 0\\ bA^{-1} & 1
\end{array}\right)
\left(\begin{array}{c} \beta \\ \chi \end{array}\right)
 +
 \left(\begin{array}{c} A^{-1}c \\ bA^{-1}c-d \end{array}\right).
$$
This verifies that the Legendre transform naturally transforms
under the action of $G$.  Moreover, the Legendre transform
coordinates $\tilde \alpha$ are independent of the potential
function $\phi$.

\subsection{Simon-Wang's developing map} \label{wang-corr}

U.\ Simon and C.P.\ Wang \cite{simon-wang93} formulate the
condition for a two-dimensional surface to be an affine sphere in
terms of the conformal geometry given by the affine metric. Since
we rely heavily on this work, we give a version of the arguments
here for the reader's convenience.  For basic background on affine
differential geometry, see Calabi \cite{calabi72}, Cheng-Yau
\cite{cheng-yau86} and Nomizu-Sasaki \cite{nomizu-sasaki}.

Consider a 2-dimensional parabolic affine sphere in $\re^3$.  Then
the affine metric gives a conformal structure, and we choose a
local conformal coordinate $z=x+iy$ on the hypersurface. Then the
affine metric is given by $h=e^{\psi}|dz|^2$ for some function
$\psi$. Parametrize the surface by $f:\mathcal{D} \rightarrow
\mathbb{R}^3$, with $\mathcal{D}$ a domain in $\mathbb{C}$. Since
$\{e^{-\frac{1}{2}\psi} f_x,e^{-\frac{1}{2}\psi} f_y\}$ is an
orthonormal basis for the tangent space, the affine normal $\xi$
must satisfy this volume condition (see e.g.\
\cite{nomizu-sasaki})
\begin{equation} \label{detxy}
\det(e^{-\frac{1}{2}\psi} f_x,e^{-\frac{1}{2}\psi} f_y, \xi) = 1,
\end{equation}
 which implies
\begin{equation}
\det(f_z,f_{\bar{z}},\xi)=\sfrac{1}{2}i e^\psi. \label{det-eq}
\end{equation}

Now only consider parabolic affine spheres.  In this case, the
affine normal $\xi$ is a constant vector, and we have
\begin{equation}
\left\{ \begin{array}{c}
D_X Y = \na_X Y + h(X,Y)\xi\\
D_X \xi = 0
\end{array} \right. \label{struc}
\end{equation}
Here $D$ is the canonical flat connection on $\re^3$, $\na$ is a
flat connection, and $h$ is the affine metric.

It is convenient to work with complexified tangent vectors, and we
extend $\na$, $h$ and $D$ by complex linearity.  Consider the
frame for the tangent bundle to the surface $\{ e_1 = f_z =
f_*(\frac{\partial}{\partial z}), e_{\bar 1}=f_{\bar z} =
f_*(\frac{\partial}{\partial {\bar z}}) \}$. Then we have
\begin{equation}
h(f_z,f_z)=h(f_{\bar z}, f_{\bar z})=0, \quad h(f_z,f_{\bar z}) =
\sfrac{1}{2}e^\psi.
\label{h-met}
\end{equation}
Consider $\theta$ the matrix of connection one-forms
$$\na e_j = \theta^k_j e_k, \quad j, k \in \{1,{\bar 1}\},$$
and $ {\hat \theta} $ the matrix of connection one-forms for the
Levi-Civita connection.   By (\ref{h-met})
\begin{equation}
{\hat \theta}^1_{\bar 1} ={\hat \theta}^{\bar 1}_1 = 0, \quad
{\hat \theta}^1_1 = \partial \psi, \quad
{\hat \theta}^{\bar 1}_{\bar 1} = {\bar \partial} \psi. \label{levi-cit}
\end{equation}

The difference ${\hat \theta} - \theta$ is given by the Pick form.  We
have
$${\hat \theta}^j_\ell - \theta^j_\ell = C^j_{\ell k} \rho^k,$$
where $\{ \rho^1 = dz,\rho^{\bar 1} = d{\bar z} \}$ is the dual
frame of one-forms.   Now we differentiate (\ref{det-eq}) and use
the structure equations (\ref{struc}) to conclude
$$\theta^1_1 + \theta^{\bar 1}_{\bar 1} = d \psi.$$
This implies, together with (\ref{levi-cit}), the apolarity condition
$$C^1_{1k} + C^{\bar 1}_{{\bar 1} k} = 0, \quad k \in \{1,{\bar 1} \}.$$
Then, when we lower the indices, the expression for the metric
(\ref{h-met}) implies that
$$C_{{\bar 1}1k} + C_{1{\bar 1}k} = 0.$$
Now $C_{jk\ell}$ is totally symmetric on three indices
\cite{cheng-yau86,nomizu-sasaki}. Therefore, the previous equation
implies that all the components of $C$ must vanish except
$C_{111}$ and $C_{ {\bar 1}{\bar 1}{\bar 1}} =
\overline{C_{111}}$.

This discussion completely determines $\theta$:
\begin{equation}
\left( \begin{array}{cc}  \theta^1_1 & \theta^1_{\bar 1} \\[1mm]
                                \theta^{\bar 1}_1 & \theta^{\bar 1}_{\bar1}
                  \end{array} \right)
 = \left( \begin{array}{cc}  \partial \psi & C^1_{{\bar 1}{\bar 1}}
                                d{\bar z} \\[1mm]
                                C^{\bar 1}_{11} dz & \bar{\partial} \psi
                  \end{array} \right)
=\left( \begin{array}{cc}  \partial \psi & \bar{U} e^{-\psi} d{\bar z} \\
                        U e^{-\psi} dz & \bar{\partial} \psi
                  \end{array} \right),
\label{conn-eq}
\end{equation}
where we define $U = C^{\bar 1}_{11} e^\psi$.

Recall that  $D$
is the canonical flat connection induced from ${\mathbb R}^3$.  (Thus,
for example, $D_{f_z}f_z = D_{\frac{\partial}{\partial z}} f_z =
f_{zz}$.)
Using this statement, together with (\ref{h-met}) and (\ref{conn-eq}),
the structure equations (\ref{struc}) become
\begin{equation}
\left\{ \begin{array}{c}
f_{zz} = \psi_z f_z + U e^{-\psi} f_{\bar z} \\
f_{{\bar z}{\bar z}} = {\bar U} e^{-\psi} f_z + \psi_{\bar z} f_{\bar z}
\\
f_{z{\bar z}} = \frac{1}{2}e^\psi \xi \end{array} \right.
\label{fzz-eq}
\end{equation}
Then, together with the equations $\xi_z=\xi_{\bar z}=0$, these
form a linear first-order system of PDEs in $\xi$, $f_z$ and
$f_{\bar z}$:
\begin{eqnarray} \label{z-deriv}
\frac{\partial}{\partial z} \left( \begin{array}{c}
\xi \\ f_z \\
f_{\bar z} \end{array} \right) &=& \left( \begin{array}{ccc}
0&0&0 \\
0& \psi_z & U e^{-\psi} \\
\frac12 e^\psi & 0 & 0
\end{array} \right)
\left( \begin{array}{c} \xi \\ f_z \\ f_{\bar z} \end{array}
\right),
\\ \label{zbar-deriv}
\frac{\partial}{\partial \bar z} \left( \begin{array}{c}
\xi
\\f_z
\\ f_{\bar z} \end{array} \right) &=& \left( \begin{array}{ccc}
0&0&0 \\
\frac12 e^\psi & 0 & 0 \\
0 & \bar{U} e^{-\psi} & \psi_{\bar z} \\
\end{array} \right)
\left( \begin{array}{c} \xi \\ f_z \\ f_{\bar z} \end{array}
\right).
\end{eqnarray}
In order to have a solution of the system (\ref{fzz-eq}),
the only condition is
that the mixed partials must commute (by the Frobenius theorem). Thus we
require
\begin{eqnarray}
\psi_{z {\bar z}} + |U|^2 e^{-2\psi}   &=& 0,
\label{psi-eq} \\
\nonumber
U_{\bar z} &=& 0.
\end{eqnarray}

The system (\ref{fzz-eq}) is an initial-value problem, in that
given (A) a base point $z_0$, (B) initial values $f(z_0)\in\re^3$,
$f_z(z_0)$ and $f_{\bar z}(z_0)=\overline{f_z(z_0)}$, and (C) $U$
holomorphic and $\psi$ which satisfy (\ref{psi-eq}), we have a
unique solution $f$ of (\ref{fzz-eq}) as long as the domain of
definition $\mathcal{D}$ is simply connected.  We then have that
the immersion $f$ satisfies the structure equations (\ref{struc}).
In order for $\xi$ to be the affine normal of $f(\mathcal{D})$, we
must also have the volume condition (\ref{det-eq}), i.e.\ $\det
(f_z, f_{\bar z}, \xi)=\frac{1}{2}ie^\psi$.  We require this at
the base point $z_0$ of course:
\begin{equation}
\det(f_z(z_0), f_{\bar z}(z_0), \xi) = \sfrac{1}{2}ie^{\psi(z_0)}.
\label{init-val}
\end{equation}
Then use (\ref{fzz-eq}) to show that the derivatives with respect
to $z$ and ${\bar z}$ of $\det(f_z,f_{\bar z},\xi)e^{-\psi}$ must
vanish. Therefore the volume condition is satisfied everywhere,
and $f(\mathcal{D})$ is a parabolic affine sphere with affine
normal $\xi$.

Using (\ref{fzz-eq}), we compute
\begin{equation} \label{U-f}
\det(f_z,f_{zz},\xi)=\sfrac12iU,
\end{equation}
 which implies that $U$
transforms as a section of $K^3$, and $U_{\bar z}=0$ means it is
holomorphic. [Note that equation (\ref{psi-eq}) is in local coordinates.
In other words, if we choose a local conformal coordinate $z$, then
the Pick form $U=U\,dz^3$, and the metric is $h=e^\psi|dz|^2$.
Then plug $U,$ $\psi$ into (\ref{psi-eq}).]

\subsection{Blaschke's holomorphic representation}
\label{blaschke-rep}

It was known to Blaschke that two-dimensional parabolic affine
spheres may be represented by two holomorphic functions (much as
minimal surfaces can).  Yau and Zaslow have recently related this
representation to the stringy cosmic string model
\cite{greene-svy90} and Hitchin's work on the moduli space of
special Lagrangian submanifolds \cite{hitchin97}. We mention two
fairly recent generalizations of this: Calabi found that affine
maximal surfaces in $\re^3$ may be represented by holomorphic data
\cite{calabi90} (an affine maximal surface is one whose area with
respect to the affine metric is at a critical point; parabolic
affine spheres are prominent examples). Also Cort\'es has
described special K\"ahler metrics in higher dimensions by
holomorphic functions \cite{cortes02} (each special K\"ahler
metric locally lives on a parabolic affine sphere of real
dimension $2n$).

We recall the version of Blaschke's result in
Ferrer-Mart\'{\i}nez-Mil\'an \cite{ferrer-mm99}. A parabolic
affine sphere with affine normal $\xi=(0,0,1)$ is given locally by
the graph $\{(\alpha^1,\alpha^2,\phi(\alpha^1,\alpha^2))\}$, where
$\phi$ is a convex function satisfying the Monge-Amp\`ere equation
$\det \frac{\partial^2\phi}{\partial \alpha^j \partial \alpha^k}
=1$. Define
\begin{eqnarray} \label{G-def}
G &=& \left(\alpha^1 + \frac{\partial\phi}{\partial
\alpha^1}\right) + i\left(\alpha^2 +
\frac{\partial\phi}{\partial \alpha^2}\right),\\
\label{F-def}
 F&=& \left(\alpha^1 - \frac{\partial\phi}{\partial
\alpha^1}\right) + i\left(-\alpha^2 + \frac{\partial\phi}{\partial
\alpha^2}\right).
\end{eqnarray}
$F$ and $G$ are holomorphic with respect to the local conformal
coordinate $z$ introduced above. Moreover the affine metric
\begin{equation} \label{holo-met}
e^\psi |dz|^2 = \sfrac14(|dG|^2 - |dF|^2)
\end{equation}
(so note that $dG\neq0$ and $|dG|>|dF|$ everywhere).

We will also need a formula for the cubic form $U$ in terms of $F$
and $G$.  Note that $f=(\alpha^1,\alpha^2,\phi)$, $\xi=(0,0,1)$.
Use (\ref{G-def}-\ref{F-def}) to write $\alpha^1$ and $\alpha^2$
in terms of $F$ and $G$.  Then compute using (\ref{det-eq})
\begin{equation} \label{U-G-F}
U= -2i \det(f_z,f_{zz},\xi) = \sfrac14(G_z F_{zz} - F_z G_{zz}).
\end{equation}

\section{Solving  \c{T}i\c{t}eica's equation} \label{solve-u}

Let $U$ be a meromorphic section of $K^3$ over $\cp^1$ with poles
of order 1 at points $p_j\in \cp^1$, and no other poles.  Assume
$U$ is not identically zero.  It is easy to see that the number of
poles is at least 6. Then we want a metric $e^\psi |dz|^2$ so that
$$ \psi_{z {\bar z}} + |U|^2 e^{-2\psi} =0,$$
in local coordinates on $M=\cp^1\setminus \{p_j\}$.  If we choose
a background metric $h$ on $M$ and then solve for a conformal
factor $e^u$ so that $e^\psi|dz|^2=e^uh$, the equation becomes
\begin{equation} \label{lap}
L_h(u)=\Delta_h u +4 e^{-2u} \|U\|_h^2 - 2\kappa_h=0.
\end{equation}
Here $\Delta_h$ is the Laplacian, $\|\cdot \|_h$ denotes the
metric on $K^3$, and $\kappa_h$ is the Gauss curvature.  Note that
for another conformal metric $k=e^v h,$
\begin{equation} \label{change-conf}
L_k(u)=e^{-v}L_h(u+v). \end{equation}

Near each pole $p_j$, we can always choose a holomorphic
coordinate $z=z_j$ so that $p_j=\{z=0\}$ and $U=\frac1z \, dz^3$
near $z=0$. For each $p_j$, we call this coordinate the
\emph{canonical holomorphic coordinate}. Then $\psi=
\log|\log|z|^2|$ solves equation (\ref{psi-eq}).  This will be our
local model near $z=0$.

\subsection{Barriers} \label{barriers}
We will solve equation (\ref{lap}) on $M$ by construction upper
and lower barriers.  We start with the lower barrier $s$, which is
easier to construct.  Let the background metric be a smooth
conformal metric $h$ so that
 \begin{equation} \label{def-h}
 h=|\log|z_j|^2||dz_j^2|
 \end{equation}
for the canonical holomorphic coordinate $z_j$ in a neighborhood
of each singularity $p_j$.

First, conformally modify $h$ so that $\tilde h = e^v h$ has
negative curvature in a neighborhood of each zero of $U$ and $v$
is compactly supported in $M$ (i.e.\ supported away from the
$p_j$).

Near each $p_j$, consider $z=z_j$, $r=|z|$, $\alpha<0$,  and
$$u=\beta|\log r|^\alpha.$$
Compute near $z=0$
\begin{equation} \label{op-u} L_h(u) = \frac1{2r^2|\log
r|^3} [ \alpha(\alpha-1)u + e^{-2u} -1 ]. \end{equation}
 For $\beta<0$, $\alpha\in(-1,0)$, then $u<0$ and it is easy to check
$L_h(u)\ge0$.

Choose $\alpha\in(-1,0)$. Let $f$ be a smooth positive function on
$M$ that is equal to $|\log |z_j||^\alpha$ on a neighborhood of
each $p_j$, and is constant outside a larger neighborhood of each
$p_j$.  In particular, we may choose $f$ to be constant in a
neighborhood of each zero of $U$.  Then if $s=v+\beta f$ for
$\beta\ll 0$, then $L_h(s)>0$ on all of $M$.

The upper barrier is harder to obtain.  Equation (\ref{op-u})
above shows that for the same model $L_h(u)\not\le0$ for
$\beta\gg0$.  In particular, the upper barrier we obtain is only
bounded at the $p_j$ (it does not go to 0).  As we'll see below in
Corollary \ref{to-zero}, the solution $u$ to $L_h(u)=0$ does go to
zero at each $p_j$.

Consider a smooth conformal background metric $k$ on $\cp^1$,
which is equal to $|dz_j|^2$ near each singularity $p_j$. Let
$\kappa_k$ be the Gaussian curvature of $k$.  Let $\tilde \kappa$
be a smooth positive function on $M$ which is equal to
$$-\sfrac12\Delta_k \log|\log|z_j|^2| = \frac1{4|z_j|^2(\log |z_j|)^2}$$
on a neighborhood of each singularity $p_j$. Note $\tilde \kappa$
is integrable. We also require that
$$\int_M (\kappa_k-\tilde \kappa) \, dV_k = 0,$$
where $dV_k$ is the volume element of the metric $k$. (This is
possible by Gauss-Bonnet.)  Now use the Green's function of the
Laplacian $\Delta_k$ to find $f$ so that
$$\Delta_k f = 2\kappa_k -2 \tilde \kappa$$
A fairly straightforward computation with the Green's function
(which we put in Appendix \ref{green-appendix} below) shows that
near each $p_j$,
$$f=\log|\log|z_j|^2| + O(1).$$
Then compute for any constant $c$
\begin{eqnarray*}
L_k(f+c) &=& \Delta_k f + 4 e^{-2f-2c}\|U\|_k^2 - 2\kappa_k \\
&=&4 e^{-2f-2c}\|U\|_k^2 - 2 \tilde\kappa.
\end{eqnarray*}
By our choice of $\tilde\kappa$, we have that if $c\gg0$,
$L_k(f+c) <0$ always.  If $h$ is the model metric from above, and
$k=e^v h$, then by (\ref{change-conf}), $L_h(f+c+v)<0$.  Let the
upper barrier $S=f+c+v$.  By the asymptotics of $f$, we know that
$S$ is bounded.  Moreover, by choosing $c$ large enough, we can
ensure that $S>s$ on $M$.

Now we solve the equation by exhausting $M$ by $M=\bigcup_n M_n$,
where $M_n$ is a smooth manifold with boundary.  For example, for
$n$ large, take
$$M_n =M \setminus \left(\bigcup_j \left\{|z_j|<\frac1n\right\} \right).$$
Let $u_n$ be the solution to the Dirichlet problem on $M_n$
$$ L_h(u)=0, \qquad u=S \quad\mbox{on}\quad \partial M_n.$$
This can be done by e.g.\ Theorem 12.5 in
\cite{gilbarg-trudinger}.  Then the strong maximum principle shows
that $S\ge u_n \ge s$, and that the sequence $\{u_n\}$ is
decreasing pointwise to a function $u$ on $M$. Therefore we have
$C^0$ bounds on $u_n$. Then equation (\ref{lap}) shows that we
have local $W^{2,p}$ bounds on $u_n$, which imply $C^{0,\alpha}$
bounds. Then (\ref{lap}) again gives local $C^{2,\alpha}$ bounds
on $u_n$ independent of $n$. Therefore, by Ascoli-Arzel\'a, the
convergence is $C^2$, and $u$ is a bounded solution of $L_h(u)=0$
on $M$.  $u$ is smooth by further bootstrapping. We record this as
a proposition.

\begin{prop} \label{u-bound}
Let $h$ be the background metric above on $M$. Then there exists a
bounded smooth solution $u$ to
$$\Delta_h u + 4 e^{-2u} \|U\|^2_h
-2 \kappa_h=0.$$
\end{prop}

\begin{rem}
Here is an interpretation of the semi-flat Calabi-Yau metric
$m=e^uh$ due to Robert Bryant.  Consider the tensor
$$ g = \frac{|U|^2}{m^2}.$$
Then $m$ is a semi-flat Calabi-Yau metric with cubic differential
$U$ if and only if $g$ has constant Gaussian curvature $-4$.  So
away from the zeroes of $U$, $g$ is a constant curvature
Riemannian metric. In the particular case where $U$ has exactly 6
poles of order 1 on $\cp^1$, then $U$ has no zeroes, and $g$ is
the unique complete metric of constant curvature $-4$ on
$\cp^1\setminus\{p_j\}$.
\end{rem}

\subsection{Blow-up analysis}
Let $z=z_j$ be the canonical holomorphic coordinate for a singular
point $p_j$.  Then for $z$ near 0,
$$h=|\log|z|^2||dz^2|,\qquad U=\frac1z\,dz^3.$$
Then $u$ satisfies
$$u_{z {\bar z}} +\frac{e^{-2u}-1}{|z|^2(\log|z|^2)^2} = 0.$$
For $\lambda\ge1$, let
$$u_\lambda(z)=u\left(\frac{z}\lambda\right).$$
 Then compute
  \begin{equation} \label{u-lambda}
 u_{\lambda,z{\bar z}} +
\frac{e^{-2u_\lambda}-1}{|z|^2(\log|z|^2-2\log \lambda)^2} = 0.
\end{equation}
 Since $u$ is bounded, this implies that $\lim_{\lambda\to\infty}
 u_{\lambda,z{\bar z}}=0$
uniformly on compact subsets of $\co \setminus \{0\}$. By the same
bootstrapping estimates as above, we have uniform $C^{2,\alpha}$
bounds on $u_\lambda$ on compact subsets of $\co \setminus \{0\}$.
Then by Ascoli-Arzel\'a and a diagonalization argument, there is a
sequence $\lambda_j\to\infty$ so that $u_{\lambda_j}$ converges
locally in $C^2$ to a limit $u_\infty$ on $\co \setminus \{0\}$.
By letting $\lambda\to\infty$ in equation (\ref{u-lambda}), we see
$$u_{\infty,z{\bar z}}=0 \quad \mbox{on}\quad \co\setminus\{0\}.$$
Therefore, $u_\infty$, as a bounded harmonic function on $\co
\setminus \{0\}$, must be a constant function.  Moreover,
$u_{\infty,z}=0$ so that $\lim_j u_{\lambda_j,z}=0$.  For any
sequence $\lambda_j\to\infty$, this argument shows that there is a
subsequence $\lambda_{j_k}$ so that $\lim_k
u_{\lambda_{j_k},z}=0$. Thus
$$\lim_{\lambda\to\infty}u_{\lambda,z}=0$$
uniformly on compact subsets of $\co \setminus \{0\}$. (We cannot
yet conclude that $u_\lambda$ converges, since the constant limit
functions for each such subsequence may be different; see
Corollary \ref{to-zero} below.) In terms of the unrescaled
solution $u$, this is equivalent to
\begin{prop} \label{u-z-limit}
$\D \lim_{z\to0}z\,u_z =0.$
\end{prop}

\section{Affine flat structure} \label{affine-flat}

Choose coordinates on $\re^3$ so that the affine normal
$\xi=(0,0,1)$. Let $\alpha^1, \alpha^2$ be coordinates on $\re^3$
transverse to $\xi$. Then locally an immersed parabolic affine
sphere is the graph of a convex function $\phi(\alpha^1,\alpha^2)$
satisfying
$$\det \frac{\partial^2\phi}{\partial \alpha^j\partial
\alpha^k}=1.$$ The affine sphere is given in coordinates by
$(\alpha^1, \alpha^2,\phi)$.

The solution of \c{T}i\c{t}eica's equation in the last section
gives induces an immersion  $f\!:\tilde M \to \re^3$, where
$\tilde M$ is the universal cover of $M$. More specifically, let
$p \in \tilde M$ with a local coordinate $z$ near $p$, and
consider an initial vector $\mathcal V \in \re^3\otimes\co$ which
satisfies
$$\det (\mathcal V,\overline{\mathcal V},\xi)=\sfrac12 i e^{\psi(p)}$$
as in equation (\ref{det-eq}) above.  Then there is a unique map
$f\!:\tilde M \to \re^3$ so that $f_z(p)=\mathcal V$ and $f$
satisfies the evolution equations
(\ref{z-deriv}-\ref{zbar-deriv}).

In the coordinates above, $f=(\alpha^1,\alpha^2,\phi)$.  If we
leave $\xi=(0,0,1)$ fixed and choose other transverse coordinates
$\check\alpha^1,\check\alpha^2$, then we have a natural affine
change of coordinates.   The map ${\rm dev}
=(\alpha^1,\alpha^2):\tilde M \to \re^2$ is the developing map of
the affine flat structure induced on $M$.

In order to check this, we must show that the holonomy
homomorphism behaves appropriately.  Let $\gamma \in \pi_1(M)$ be
a deck transformation. Consider a loop representing $\gamma$ which
begins at a point $q$, and choose a lift $\tilde q\in\tilde M$.
Lift the loop so that the other endpoint is $\tilde q'\in\tilde
M$.  Then following a neighborhood of the lifted path induces a
coordinate map from a neighborhood of $\tilde q$ to a neighborhood
of $\tilde q'$.  This map is a constant element of $\Sl{2}\ltimes
\re^2$, by the affine invariance of the initial value problem: Let
$\Phi$ be the element of $\Sl{3}\ltimes\re^3$ which takes the
initial data consisting of the position $f$ and the frame
$\{f_z,f_{\bar z},\xi\}$ at $\tilde q$ to the corresponding data
at $\tilde q'$. Then since $f$ solves the initial value problem
with initial data at $\tilde q$, $\Phi(f)$, $\Phi_*\{f_z,f_{\bar
z},\xi\}$ must solve the initial value problem with initial data
at $\tilde q'$. Thus $\Phi\circ f = f \circ \gamma$ for the deck
transformation $\gamma$. Equation (\ref{det-eq}) shows that $\Phi$
preserve volume. Since the affine normal $\xi$ is constant, $\Phi$
induces an affine map ${\rm hol_\gamma}\in \Sl{2}\ltimes\re^2$.
Thus we have a holonomy map $\gamma\mapsto {\rm hol}_\gamma$ from
$\pi_1(M)$ to $\Sl{2}\ltimes\re^2$. The pair $\{{\rm dev},{\rm
hol}\}$ is equivalent to the affine flat structure (see e.g.\
Goldman \cite{goldman88}).  (We should also show that if we choose
a different basepoint in $\tilde M$ and different coordinates in
$\re^3$ transverse to $\xi$, that the pair $\{{\rm dev},{\rm
hol}\}$ transforms appropriately.  These points are again easy to
check by the affine invariance of the problem in $\re^3$ and the
uniqueness of solutions to the initial value problem.)

We have proved

\begin{prop}
The immersion $f\!:\tilde M \to \re^3$ induces an affine flat
structure on $M$ with a covariant constant volume form.
\end{prop}

\subsection{The holonomy type}

We focus on the affine structure around each pole of the cubic
differential $U$.  Choose a loop around such a singular point,
represented by an element $\gamma$ of $\pi_1(M)$.

Near a singularity $z=z_j=0$, let $z=e^{iw}$ so that for
$|z|\in(0,\epsilon)$, $w=x+iy$ satisfies $y>-\log\epsilon$.  Then
in terms of the $w$ coordinates,
\begin{eqnarray*}
e^\psi |dz|^2 &=&e^u|\log|z|^2||dz|^2= 2ye^{-2y}e^u|dw|^2,\\
U&=&\frac1z\,dz^3 =-ie^{2iw}\,dw^3.
 \end{eqnarray*}
 Now the equations (\ref{z-deriv}), (\ref{zbar-deriv}) imply that
in terms of the real frame $\{\xi,f_x,f_y\},$
 \begin{equation} \label{conn-x}
 \left(
\begin{array}{c}
\xi \\ f_x \\
f_y \end{array} \right)_x = A\left( \begin{array}{c} \xi \\ f_x \\
f_y \end{array} \right),
\end{equation}
 where $A=A(x,y)=$
 $$ \left( \begin{array}{ccc}
0&0&0 \\
2ye^{-2y}e^u& \frac12u_x+e^{-u}\frac1{2y}\sin2x &
  -\frac12u_y-\frac1{2y}+1+e^{-u}\frac1{2y}\cos2x \\
0 & \frac12u_y+\frac1{2y}-1+e^{-u}\frac1{2y}\cos2x &
  \frac12u_x -e^{-u}\frac1{2y}\sin2x
\end{array} \right).
$$
In the $w$ coordinate, Proposition \ref{u-z-limit} implies that
$$u_w=\sfrac12u_x-\sfrac{i}2u_y\to0 \quad \mbox{as}\quad
y\to\infty.$$
 Therefore, since $u$ is bounded, we find
 $$\lim_{y\to\infty}A=A_\infty=\left(\begin{array}{ccc}
 0&0&0\\0&0&1\\0&-1&0\end{array}\right)$$
 uniformly in $x$.

To compute the affine flat structure, we only need the components
of $f$ transverse to $\xi$, which we call ${\rm
dev}=(\alpha^1,\alpha^2)$ above.  Then, we only need to consider
$\tilde A$, the bottom right $2\times2$ submatrix of $A$ and we
have
$$ \left(\begin{array}{c} {\rm dev}_x \\ {\rm dev}_y
\end{array}\right)_x = \tilde A \left(\begin{array}{c} {\rm dev}_x \\ {\rm dev}_y
\end{array}\right), \qquad \lim_{y\to\infty}\tilde A = \tilde A_\infty=
\left(\begin{array}{cc}0&1\\-1&0\end{array}\right).$$

The theory of ODEs \cite{hartman} then implies that the
fundamental solution to the ODE (\ref{conn-x}) must approach the
fundamental solution to
 \begin{equation} \label{const-coeff-x} X_x=\tilde A_\infty X,
 \end{equation}
where $X=({\rm dev}_x,{\rm dev}_y)^\perp$. In other words,  the
solution $\Psi(x,y)$ to the initial value problem
$$ X(0,y)= X_0,\qquad X_x=\tilde A X,$$
must approach the solution
$$X=X_0 e^{x \tilde A_\infty}$$
to (\ref{const-coeff-x}) uniformly in $x$ as $y\to\infty$.

For any $y\gg0$, the linear path from $(0,y)$ to $(2\pi,y)$
corresponds to a loop $|z|=e^{-y}$ around the singularity $p_j$.
Also, $\{{\rm dev}_x,{\rm dev}_y\}$ is a frame of the tangent
space to hypersurface $H$. Therefore, integrating the initial
value problem (\ref{z-deriv}-\ref{zbar-deriv}) along such path
computes the linear part of the holonomy.  So the solution to the
ODE $X_x=\tilde A X$ computes this part of the holonomy. Note that
since the deck transformation is $w\mapsto w+2\pi$, the frame
$\{{\rm dev}_x,{\rm dev}_y\}$ is appropriate for computing the
holonomy (the frames at $x=0$ and $x=2\pi$ may be naturally
identified).

Since the connection $D$ on $\re^3$ is flat, the conjugacy class
of the holonomy is independent of the choice of loop in a free
homotopy class.  In particular, the linear holonomy matrix
determined by our frame and the loop $|z|=e^{-y}$ is given by
$\Psi(2\pi,y)$, which satisfies
$$\lim_{y\to\infty}\Psi(2\pi,y)=e^{2\pi \tilde A\infty}=
\left( \begin{array}{cc} 1&0 \\ 0&1 \end{array} \right).$$ This
\emph{does not} mean the linear part of the holonomy is trivial,
as shown by the family of matrices
$$\left(\begin{array}{cc} 1 & \epsilon \\ 0&1\end{array}\right)
\to \left( \begin{array}{cc} 1&0 \\ 0&1 \end{array} \right).$$
However, we do know that the set of eigenvalues of the linear part
of the holonomy must be $\{1\}$.  We record this as

\begin{prop}
The set of eigenvalues for the linear part of the affine holonomy
around each puncture is $\{1\}$.
\end{prop}

There are four different conjugacy classes in $\Sl{2}\ltimes\re^2$
whose elements have unipotent linear part (we simply list a
representative in each conjugacy class):
\begin{enumerate}
\item $Y \mapsto Y$.

\item $Y \mapsto Y +b, \quad b\neq 0$. \label{trans}

\item $Y \mapsto
\left(\begin{array}{cc}1&1\\0&1\end{array}\right)Y$.

\item $Y \mapsto
\left(\begin{array}{cc}1&1\\0&1\end{array}\right)Y + b,\qquad$
which has no fixed point. \label{no-fixed}
\end{enumerate}

In order to address the affine part of the holonomy, we must
consider dev itself instead of the derivatives ${\rm dev}_x$,
${\rm dev}_y$.  To calculate ${\rm dev}=\int {\rm dev}_y\,dy$,
compute from the structure equations
(\ref{z-deriv}-\ref{zbar-deriv})
\begin{equation} \label{conn-y}
 \left(
\begin{array}{c}
\xi \\ f_x \\
f_y \end{array} \right)_y = B\left( \begin{array}{c} \xi \\ f_x \\
f_y \end{array} \right),
\end{equation}
 where $B=B(x,y)=$
 $$ \left( \begin{array}{ccc}
0&0&0 \\
0 & \frac12u_y+ \frac1{2y} -1 + e^{-u}\frac1{2y}\cos2x &
  \frac12u_x- e^{-u}\frac1{2y}\sin2x \\
2ye^{-2y}e^u & -\frac12u_x - e^u\frac1{2y}\sin2x &
  \frac12u_y+\frac1{2y}-1 - e^{-u}\frac1{2y}\cos2x
\end{array} \right).
$$
As above, let $\tilde B$ be the bottom right $2\times2$ submatrix
of $B$, and we have
$$ \left(\begin{array}{c} {\rm dev}_x \\ {\rm dev}_y
\end{array}\right)_y = \tilde B \left(\begin{array}{c} {\rm dev}_x \\ {\rm dev}_y
\end{array}\right), \qquad \lim_{y\to\infty}\tilde B = \tilde B_\infty=
\left(\begin{array}{cc}-1&0\\0&-1\end{array}\right).$$

A theorem of Perron \cite{perron29} and Lettenmeyer
\cite{lettenmeyer} (see also Hartman-Wintner
\cite{hartman-wintner55}) then gives the asymptotic behavior of
any initial value problem $X(y_0)=X_0$, $X_y=\tilde By$:
$$|X| = e^{-y+o(y)} \quad \mbox{as} \quad y\to\infty.$$
So then if we  choose $X=({\rm dev}_x, {\rm dev}_y)^\perp$ as above,
with initial conditions
$X(x_0,y_0)=X_0$, then as $y\to \infty$, $|X(x_0,y)|=e^{-y+o(y)}$.
This gives the behavior of integrating
(\ref{z-deriv}-\ref{zbar-deriv}) along a path $x=x_0$. If we have
followed this path from $(x_0,y_0)$ to $(x_0,y)$, we can then
integrate in the $x$ direction to find that
$$X(x,y) = \left[\left( \begin{array}{cc}
\cos x & \sin x \\ -\sin x & \cos x \end{array} \right) + o(1)
\right] X(x_0,y)$$ as $y\to\infty$ and bounded $x$.   Therefore,

\begin{lem}
The solution $X$ to the initial value problem $X(x_0,y_0)=X_0$,
$X_x=\tilde A X$, $X_y = \tilde B X$, satisfies $|X| =
e^{-y+o(y)}$ as $y\to \infty$ uniformly in any bounded interval in
$x$.
\end{lem}

Therefore $|{\rm dev}_y| = e^{-y+o(y)}$ as $y\to \infty$ uniformly
in a bounded interval of $x$ (and the same is true for ${\rm
dev}_x$). Then if we augment the initial value problem to include
an initial value for ${\rm dev}(x_0,y_0)$, then
$${\rm dev}(\infty) = {\rm dev}(x_0,y_0)+
\int_{y_0}^\infty {\rm dev}_y\, dy
$$
exists. ($x=x_0$ is constant in the path of this integral.)
Moreover, we can determine how the holonomy acts on ${\rm
dev}(\infty)$ by solving the initial value problem along the
following path:  First let $y$ go from $y_0$ to $y'$ for some
$y'\gg0$. Then let $x$ go from $x_0$ to $x_0+2\pi$, and finally
let $y\to \infty$.  By the decay of ${\rm dev}_x$, the
contribution from integrating from $x_0$ to $x_0+2\pi$ is of order
$e^{-y'+o(y')}$. Therefore, the limiting point value of dev for
such a path must be ${\rm dev}(\infty)$ again as $y'\to\infty$.
Moreover, the integral is the same upon integrating along any two
homotopic paths; thus the limiting case as $y'\to\infty$ is equal
to the integral along any such path for $y'\gg0$. These paths
compute the action of the holonomy; so ${\rm dev}(\infty)$ is a
fixed point of the affine holonomy. We record this as

\begin{prop} The affine holonomy around each puncture has a fixed
point.
\end{prop}

Note that this rules out the holonomy cases (\ref{trans}) and
(\ref{no-fixed}) above. Thus the affine holonomy either is the
identity or is conjugate to  $Y\mapsto
\left(\begin{array}{cc}1&1\\0&1\end{array}\right)Y$.

To rule out the identity holonomy, choose affine flat coordinates
$\alpha^1$, $\alpha^2$, and recall the representation in terms of
local holomorphic functions $F$ and $G$. The definitions
(\ref{G-def}-\ref{F-def}) imply that
$$\alpha^1 = \sfrac12\,{\rm Re}\,(G+F), \qquad \alpha^2 = \sfrac12
\,{\rm Im}\,(G-F).$$
 Since the holonomy is trivial, $\alpha^1$ and $\alpha^2$ are
well-defined on some neighborhood of the puncture
$\{z:0<|z|<\epsilon\}$.  Therefore,  $\frac{\partial}{\partial
z}(G+F) = 4 \frac{\partial \alpha^1}{\partial z}$ is single
valued. This is similarly true for $\frac{\partial}{\partial
z}(G-F)$, and so $dF/dz$ and $dG/dz$ are single valued. We will
derive a contradiction given the bounds on the metric. By formula
(\ref{holo-met}) and Proposition \ref{u-bound} above, the metric
$e^\psi|dz|^2$ satisfies
\begin{equation} \label{grow-met}
C'|\log|z|^2| |dz|^2\le \sfrac14(|dG|^2 - |dF|^2) \le C|\log|z|^2|
|dz|^2
\end{equation}
 for $C,C'>0$. In particular, $dG/dz$ cannot
have an essential singularity since it goes to infinity as
$z\to0$; so it must have a pole of some order $n$ at $z=0$.
$dF/dz$ then cannot have an essential singularity, since it must
satisfy $|dF/dz|^2\ge |dG/dz|^2 - 4C|\log|z|^2|$, which forces it
to have a pole. Then look at the power series of $dG/dz$ and
$dF/dz$ to derive a contradiction. Therefore, the affine holonomy
cannot be trivial and we have proved

\begin{thm}
For any oriented loop around each pole $p_j$ of $U$, the affine
holonomy corresponding to the metric constructed in Proposition
\ref{u-bound} is conjugate to
$$Y \mapsto \left(\begin{array}{cc}1&1\\0&1\end{array}\right)Y.$$
\end{thm}

\begin{rem}
We call the holonomy type in the previous theorem
\emph{parabolic}.
\end{rem}

\subsection{Fine structure}

To investigate the fine structure, we again write $dF/dz$ and
$dG/dz$ in terms of holomorphic functions, this time involving log
terms. We will use the following terminology: A function is
\emph{holomorphic on the disk} if it is holomorphic on
$\{z:|z|<\epsilon\}$ for $\epsilon$ a small positive number.
Similarly, a function is \emph{holomorphic on the punctured disk}
if it is holomorphic on $\{z:0<|z|<\epsilon\}$. The positive
constant $\epsilon$ will be unspecified, and it may be shrunk a
little if necessary.

Introduce coordinates $\{\alpha^1,\alpha^2\}$ so that the holonomy
takes
\begin{equation} \label{hol-alpha}
 \alpha^1 \mapsto \alpha^1 + \alpha^2, \qquad \alpha^2\mapsto
\alpha^2,
\end{equation}
and so that for any fixed $x$,
\begin{equation} \label{alpha-lim}
\lim_{y\to\infty}\alpha^1 =\lim_{y\to\infty}\alpha^2 = 0.
\end{equation}
Note  these last equations are equivalent to ${\rm
dev}(\infty)=0$. Then as above, since $\alpha^2$ is single-valued
and ${\rm Im}\,(G-F) = 2\alpha^2$, we find
 \begin{equation} \label{G-Fprime} (G-F)'(z) = j(z),\end{equation}
for $j(z)$ holomorphic on the punctured disk. Since ${\rm
Re}\,(G+F) = 2\alpha^1$, (\ref{hol-alpha}) and (\ref{G-Fprime})
imply that
$$ (G+F)'(z) + \frac{\log z}{2\pi} j(z)$$
is holomorphic on the punctured disk.  Therefore,
$$ G'(z) = -\sfrac1{4\pi}j(z) \log z + \ell(z), \qquad
F'(z) =  -\sfrac1{4\pi}j(z) \log z + \ell(z) - j(z).$$
 for $\ell$ holomorphic on the punctured disk.

Since the metric  is $\frac14(|dG|^2-|dF|^2)$,
\begin{equation} \label{met-j}
 |G'|^2 - |F'|^2 = -\sfrac1{4\pi}|j|^2 \log|z|^2 +
\ell\bar\jmath + \bar\ell j  - |j|^2
\end{equation}
 is always positive. In particular, $j$ can never be 0. Also notice that
$G'$ is never zero on the punctured disk, even though it's not
single-valued. Then for
$$k(z) = -\frac{4\pi\ell(z)}{j(z)},$$
$\log z + k = -4\pi G'/j$ is never zero. Exponentiating, we find
$z e^k \neq 1$ on the punctured disk.  Of course $ze^k\neq 0$
there as well. Picard's Big Theorem then implies $ze^k$ (and thus
also $e^k$) cannot have an essential singularity at $z=0$.  Now
write $k=k_1 + k_2$, for
 \begin{equation} \label{series-rep}
  k_1 = \sum_{n=1}^\infty a_{-n} z^{-n},
\qquad k_2 = \sum_{n=0}^\infty a_n z^n.
\end{equation}
 We want to show that $k_1=0$, (i.e.\ $k$ extends to a holomorphic
function on the disk).  Since  $e^{k_2}$ is nonzero and
holomorphic on the disk, $e^{k_1}=e^k/e^{k_2}$ then cannot have an
essential singularity at $z=0$, and $e^{k_1}$ must be a polynomial
in $1/z$. Then
$$e^{k_1(z)} = z^{-n} p(z),\qquad k_1(z) = -n\log z + \log p(z)$$
for $p$ a polynomial in $z$ which doesn't vanish on a neighborhood
of $z=0$. Unless $n=0$ (and thus $k_1=0$), $k_1$ cannot be written
as a Laurent series convergent on the punctured disk. This
contradicts (\ref{series-rep}); so $k_1=0$, and $k$ extends over
$z=0$ as a holomorphic function.

Equations (\ref{grow-met}) and (\ref{met-j}) then give the
following bounds on the metric $e^\psi |dz|^2$:
$$ C'|\log|z|^2| \le e^\psi =
\sfrac14|j|^2(-\sfrac1{4\pi} \log |z|^2 - \sfrac1{2\pi}\,
{\rm Re}\, k -1) \le C|\log|z|^2|
$$ Since $k$ is holomorphic on the disk, this implies $j$ is bounded near
$z=0$ and thus $j$ extends over $z=0$ as a nonvanishing
holomorphic function.  So
 \begin{equation} \label{G-F}
  G' = -\sfrac j{4\pi}(\log z + k),
\qquad F' = -\sfrac j{4\pi} (\log z + k + 4\pi).
 \end{equation}
We also know that the cubic form $U = 1/z\,dz^3$, and thus
(\ref{U-G-F}) gives
$$\frac1z = \sfrac14(G'F''-G''F') = -\frac{j^2}{16\pi} \left(\frac1z
+ k'\right).$$
 Thus we have proved

\begin{prop} \label{holo-func-prop}
Consider the canonical coordinate $z=z_j$ around each pole $p_j$
of $U$, and the semi-flat Calabi-Yau metric $e^uh$ constructed in
Proposition \ref{u-bound} above. For affine flat coordinates
$\alpha^1$, $\alpha^2$ satisfying
(\ref{hol-alpha}-\ref{alpha-lim}), the holomorphic functions $G$
and $F$ defined in (\ref{G-def}-\ref{F-def}) depend on the
following data: Let $j(z)$, $k(z)$ be two holomorphic functions on
a neighborhood of $p_j=\{z=0\}$ so that
\begin{equation} \label{j-k}
-16\pi = (1 + zk')j^2. \end{equation}
  Then $G$ and $F$ satisfy
(\ref{G-F}).
\end{prop}

\begin{rem}
The normalization (\ref{alpha-lim}) also implies that along any
radial path toward the origin $z=0$, ${\rm Im}(G-F)$ and ${\rm
Re}(G+F)$ go to zero.
\end{rem}

\begin{cor} \label{to-zero}
The conformal factor $u$ constructed in Proposition \ref{u-bound}
above approaches 0 at each of the poles of the cubic form $U$.
\end{cor}

\begin{proof}
 By (\ref{j-k}), $j^2(0) = - 16\pi$ and so the metric
$$ e^\psi = \sfrac14|j|^2(-\sfrac1{4\pi} \log |z|^2 - \sfrac1{2\pi}\,
{\rm Re}\, k -1) = |\log|z|^2| + O(1).
$$
But near $z=0$, (\ref{def-h}) implies $e^\psi = e^u|\log|z|^2|$.
\end{proof}

The description provided by Proposition \ref{holo-func-prop}
allows us to calculate another important invariant of the affine
structure near each pole of $U$: a winding number.

Consider an affine flat structure on a punctured disk with
parabolic holonomy and coordinates $(\alpha^1,\alpha^2)$
satisfying (\ref{hol-alpha}-\ref{alpha-lim}).  Then the line
$L=\{\alpha^2=0\}$ is preserved by the action of the holonomy. If
possible, choose a point $p$ in the preimage of $L$ under the
developing map.  Then choose a path $\mathcal P$ starting and
ending at $p$ which winds once around the puncture of the
punctured disk. Since ${\rm dev}(p)\in L$, the developing map
takes any lift of $\mathcal P$ to a path $\rm{dev}(\mathcal P)$ in
$\re^2\setminus \{0\}$ which begins and ends at same point. We
define the \emph{winding number} of the affine flat structure to
be the winding number of ${\rm dev}(\mathcal P)$ around the
origin. If the developing map of the punctured disk does not
intersect $L$, we define the winding number to be zero.

\begin{thm}
For the affine flat structure constructed above, the winding
number around each pole of $U$ is $+1$.
\end{thm}

\begin{proof}
Proposition \ref{holo-func-prop} gives us that
\begin{eqnarray*}
\alpha^1 &=&\sfrac12 \,{\rm Re}\, (G+F) \\
 &=& \sfrac12 \,{\rm Re}\, \int -\sfrac{j}{2\pi} ( \log z + k +
      2\pi)\,dz \\
 &=& -\sfrac1{4\pi} \,{\rm Re} \,\int ic\log z + a + O(z\log z)\,dz \\
 &=& \frac{cr}{4\pi} [ (\theta + b_1) \cos\theta + (\log r -
 b_2)\sin\theta ] + O(z^2\log z).
\end{eqnarray*}
Here $z=re^{i\theta}$, $c = j(0)=\pm 4\sqrt\pi$, $a=(k(0)+2\pi)c$,
and $b_1+ib_2 = ic+a$. There is no constant term in the
integration by assumption (\ref{alpha-lim}). Similarly,
$$ \alpha^2 = \frac{cr}2 \cos \theta + O(z^2).$$

Fix $r$ near 0.  Then, as we show below, we may ignore higher
order terms, and the winding number is easily computed to be $+1$
for
$$\check\alpha = (\check\alpha^1,\check\alpha^2)=\left(\frac{cr}{4\pi}
(\log r - b_2)\sin \theta, \frac{cr}2\cos
\theta\right).$$
 More specifically, assume $|\log r -b_2| \gg 2\pi$. Then choose
$\Theta \in \frac\pi2 + 2\pi \mathbb Z$ so that $\Theta + b_1 \ll
|\log r - b_2|$.  Then for $\theta \in [\Theta,\Theta+2\pi]$, the
path $\check\alpha(\theta)$ winds around the origin once with
orientation (note $\log r - b_2<0$), and
$\check\alpha(\Theta)=\check\alpha(\Theta+2\pi)$ is in the fixed
line $L$.

Since as $r\to0$, $\alpha=(\alpha^1,\alpha^2)$ is $C^1$ close to
$\check\alpha$, there is a $\delta$ near zero so that
$\alpha(\Theta + \delta)$ is in the fixed line $L$, and the curve
$\alpha(\theta)$ meets $L$ transversely at $\theta=\Theta+\delta$.
Moreover, the winding number of $\alpha$ for
$\theta\in[\Theta+\delta,\Theta+\delta+2\pi]$ is $+1$, since the
loops determined by $\alpha$ and $\check\alpha$ are homotopic via
$$ \alpha_t(\theta)=t\check\alpha(\theta) + (1-t)\alpha(\theta-\delta), \qquad
\theta\in[\Theta,\Theta+2\pi].$$
 For ranges of $\theta\in[\Theta+\delta,\Theta+\delta+2\pi] +
2\pi m$, $m\in\mathbb Z$, the winding number is still $+1$ since
the path $\alpha(\theta)$ is just shifted $m$ times by the action
of the holonomy.
\end{proof}

All together, we have the following description of the metric and
affine flat structure.

\begin{thm} \label{big-theorem}
Given a holomorphic cubic differential $U$ on $\cp^1$ with poles
of order 1 at $\{p_j\}$, there is a semi-flat Calabi-Yau metric on
$\cp^1\setminus\{p_j\}$ asymptotically given by
$$\left[|\log|z_j|^2| + o(1)\right]\,|dz_j|^2$$
for $z_j$ the canonical holomorphic coordinate near the pole
$p_j$.  The affine flat structure has parabolic holonomy and
winding number $+1$ around $p_j$.
\end{thm}

\begin{rem}
In terms of the affine flat coordinates near the singularity
$p_j$, the asymptotics of the metric and the affine K\"ahler
potential function $\phi$ can also be easily calculated from the
holomorphic representation in Proposition \ref{holo-func-prop}.
\end{rem}

\section{Mirror symmetry} \label{mirror}

Consider the picture of Strominger-Yau-Zaslow in the simplest case
(without instanton corrections or singular fibers).  A Calabi-Yau
manifold $X$ admits a map
 \begin{equation} \label{fibration}
 \pi\: X \to \bar B,
 \end{equation}
where $B$ is an affine flat manifold and $\pi$ is a fibration with
fiber $T^n$ of special Lagrangian $n$-tori. $B$ admits a semi-flat
Calabi-Yau metric.  Over an affine coordinate chart $\Omega\subset
B$, form the tube domain $\Omega + i\re^n$. The special Lagrangian
tori fibered over $\Omega$ are then quotients of the imaginary
fibers $i\re^n$.  The mirror Calabi-Yau manifold should then be
constructed by taking a Fourier transform in the fiber variables
and a Legendre transform in affine coordinates and affine K\"ahler
potential on $B$.  The details of this construction may be found
in Leung \cite{leung00,leung02}.

The semi-flat Calabi-Yau metric we construct is singular at the
poles $p_j$ of $U$.  As Gross-Wilson showed \cite{gross-wilson00},
a semi-flat metric with similar behavior at the singularities is
obtained as the Gromov-Hausdorff limit of certain classes of
elliptic K3 surfaces equipped with Calabi-Yau metrics. Near the 24
singular points, their model for nearby smooth Calabi-Yau metrics
is not semi-flat (they glue in a model metric due to Ooguri-Vafa
\cite{ooguri-vafa96} there).  In particular, the fibration
(\ref{fibration}) is not globally valid for any Calabi-Yau
manifold near this singular limit.  Therefore, the Fourier
transform on the fibers is not relevant in our case, and mirror
symmetry expresses itself only through the Legendre transform of
the affine coordinates and the semi-flat Calabi-Yau potential
function.

The Legendre transform appears naturally in the holomorphic
representation in Subsection \ref{blaschke-rep} above.  In
particular the dual affine coordinates $\beta_j$ under the
Legendre transform are by (\ref{G-def}-\ref{F-def})
$$ \beta_1=\frac{\partial \phi}{\partial\alpha^1} = \sfrac12{\rm
Re}\,(G-F), \qquad \beta_2 = \frac{\partial\phi}{\partial
\alpha^2} = \sfrac12 {\rm Im}\,(G+F).$$
 The Legendre transform potential is given by
 $$\chi = \alpha^1\beta_1 + \alpha^2\beta_2 - \phi =
 \sfrac14(|G|^2-|F|^2) - \phi.$$
In terms $\alpha^j$ and $\beta_j$, $G$ and $F$ are given by
$$ G = (\alpha^1+\beta_1) + i(\alpha^2 + \beta_2), \qquad
F = (\alpha^1 - \beta_1) + i(\alpha^2 - \beta_2).$$
 Passing a semi-flat Calabi-Yau metric to its mirror means that we
switch the roles of the $\alpha^j$ and the $\beta_j$.  So the
mirror transform becomes
$$ G\mapsto G,\qquad F\mapsto -F.$$
The metric remains the same under this mapping, and the cubic
differential transforms as
$$ U \mapsto -U$$
by (\ref{U-G-F}).

We can also use Proposition \ref{holo-func-prop} to find explicit
asymptotics of the affine flat coordinates $\beta_j$ and potential
$\chi$.  In particular, similarly to (\ref{alpha-lim}), we may
assume that
 \begin{equation} \label{beta-to-zero}
\beta_1,\beta_2\to0 \end{equation}
 along any radial path as
$z\to0$ for $z$ the canonical holomorphic coordinate along each
pole. Since $\beta_j=\partial\phi/\partial\alpha^j$, this may be
accomplished by finding a new tilted set of coordinates in $\re^3$
of the type allowed in Subsection \ref{aff-flat-coord} above. So
together the assumptions (\ref{alpha-lim}) and
(\ref{beta-to-zero})  are equivalent to requiring the holomorphic
functions $F$ and $G\to0$ along any radial path as $z\to0$ (see
the remark after Proposition \ref{holo-func-prop}). Moreover, we
can read off the conjugacy class of the holonomy and the winding
number to conclude

\begin{prop} The affine flat structure for the
mirror semi-flat Calabi-Yau has parabolic holonomy and winding
number $+1$ around each pole of the cubic differential $U$.  It is
naturally isometric to its mirror.
\end{prop}

\appendix
\section{Green's function calculation} \label{green-appendix}

Let
$$f(p) = \int_{S^2} G(p,q)
\,[2\kappa_k(q)-2\tilde\kappa(q)]\,dV_k(q)$$ for $k$, $\kappa_k$,
$\tilde \kappa$ defined as in subsection \ref{barriers} above and
$G(p,q)$ the Green's function with respect to $\Delta_k$.  Then in
the coordinate $z_j$ near each pole $p_j$ of $U$,

\begin{lem}
$f(z_j) = \log|\log|z_j|^2| + O(1)$.
\end{lem}

\begin{proof}
The Green's function on a compact Riemannian surface is of the
form
$$G(P,Q) = \frac1{2\pi} \log d(P,Q) + O(1)$$
for $d$ the Riemannian distance. Thus
$$f(P) = \int_{S^2} G(P,Q)
\,[2\kappa_k(Q)-2\tilde\kappa(Q)]\,dV_k(Q).$$
 Represent $P$ by the coordinate $z_j$ near $z_j=0$ only.  Then
for a small positive $\delta$, $P= re^{i\theta}$, and $Q = \rho
e^{i\varphi}$,
\begin{eqnarray*}
f(P) &=& O(1) + g(r), \\
g(r) &=& \int_{|Q|<\delta}\left(\frac1{2\pi}\log|P-Q|\right)
 \left( -\frac1{2|Q|^2(\log|Q|)^2} \right)\,dq_1dq_2\\
&=&  - \frac1{4\pi} \int_{|Q|<\delta} \frac{\log|re^{i\theta}-\rho
e^{i\varphi}|}{\rho^2 (\log\rho)^2} \, \rho \,d\rho d\varphi.
\end{eqnarray*}
We may change variables of integration to assume $\theta=1$ and
compute
\begin{eqnarray*}
g'(r) &=&  - \frac1{4\pi}
 \int_{|Q|<\delta} \frac{\frac{\partial}{\partial r}
 \log \sqrt{(r-\rho\cos\varphi)^2 + \rho^2\sin^2\varphi}} {\rho
 (\log\rho)^2}\,d\rho d\varphi \\
&=& - \frac1{4\pi} \int_0^\delta \frac1{\rho(\log\rho)^2}
    \int_0^{2\pi} \frac{2r -2 \rho\cos\varphi}{r^2+\rho^2 -2r\rho
    \cos\varphi}\, d\varphi \,
    d\rho.
\end{eqnarray*}

The inner integral can be integrated \cite[\S
2.554]{gradshteyn-ryzhik}, but we may also evaluate it as a
contour integral for a new complex variable $\zeta = \chi
e^{i\varphi}$: For $|\zeta|=1$, $\cos \varphi = \frac12(\zeta +
\frac1\zeta)$ and thus
$$
 \int_0^{2\pi} \frac{2r -2 \rho\cos\varphi}{r^2+\rho^2 -2r\rho
    \cos\varphi}\, d\varphi = \int_{|\zeta|=1} \frac{\rho
    \zeta^2 - 2r\zeta + \rho} {r\rho \zeta^2 -(r^2+\rho^2) \zeta +
    r\rho} \, \left(-i\frac{d\zeta}\zeta\right) $$
 There are poles of the integrand at $\zeta = 0,$ $\rho/r$ and
 $r/\rho$ with residues $-i/r$, $-i/r$ and $i/r$ respectively.
 Thus if $\rho<r$, the sum of the residues inside the contour is
 $-2i/r$, and if $r<\rho$, the sum of the residues inside the
 contour is 0. Therefore,
$$
\int_0^{2\pi} \frac{2r -2 \rho\cos\varphi}{r^2+\rho^2 -2r\rho
    \cos\varphi}\, d\varphi = \left\{
    \begin{array}{c@{\quad \mbox{if}\quad}l}
      \frac{4\pi}{r} & \rho<r \\
      0 & r<\rho,
    \end{array} \right.
$$
and for $r<\delta$,
$$
g'(r) = -\frac1r \int_0^r \frac{d\rho}{\rho (\log \rho)^2} =
\frac1{r\log r}.
$$
 Therefore, $$f(P) = O(1) + g(r) = O(1) + \log|\log r| =
 O(1) + \log|\log |P|^2|.$$
\end{proof}

\bibliographystyle{abbrv}
\bibliography{thesis}

\end{document}